\documentclass{amsart}

\title{Spectral properties of bipolar minimal surfaces in $\mathbb S^4$}
\author{Hugues Lapointe}
\date{\today}

\numberwithin{equation}{subsection}
\theoremstyle{definition}
\newtheorem{definition}[equation]{Definition}
\theoremstyle{plain}
\newtheorem{lemma}[equation]{Lemma}
\newtheorem{theorem}[equation]{Theorem}

\newtheorem{prop}[equation]{Proposition}

\begin{document}

\begin{abstract}
The $i$-th eigenvalue of the Laplacian on a surface can be viewed as a functional on
the space of Riemannian metrics of fixed area.
Extremal points of these functionals correspond to surfaces admitting minimal isometric
immersions into spheres. Recently, critical metrics for the first eigenvalue were classified
on tori and on Klein bottles. The present paper is concerned with  extremal metrics for higher
eigenvalues on these surfaces. We apply a classical construction due to Lawson. For the bipolar
surface $\tilde \tau_{r,k}$ of the Lawson's torus or Klein bottle $\tau_{r,k}$ it is shown that:
\begin{enumerate}
\item If $rk \equiv 0 \mod 2$, $\tilde \tau_{r,k}$ is a torus with an extremal metric for $\lambda_{4r - 2}$.
\item If $rk \equiv 1 \mod 4$, $\tilde \tau_{r,k}$ is a torus with an extremal metric for $\lambda_{2r - 2}$.
\item If $rk \equiv 3 \mod 4$, $\tilde \tau_{r,k}$ is a Klein bottle with an extremal metric for $\lambda_{r - 2}$.
\end{enumerate}
Furthermore, we find explicitly the $\mathbb S^1$-equivariant minimal immersion of the bipolar surfaces
into $\mathbb S^4$ by the corresponding eigenfunctions.
\end{abstract}

\maketitle


\section{Introduction and main results}

\subsection{Extremal metrics for eigenvalues}

Let $\Delta$ be the Laplacian on a closed surface $M$ with metric $g$.  In local coordinates, if $g = \Sigma g_{ij} dx_i dx_j$,
\[
\Delta_g f = - \frac{1}{\sqrt{|g|}} \Sigma \frac{\partial}{\partial x^i} \left( \sqrt{|g|} g^{ij} \frac{\partial f}{\partial x^j} \right)
\]
Let
\[
0 = \lambda_0(M,g) < \lambda_1(M,g) \leq \lambda_2(M,g) \leq \lambda_3(M,g) \leq \hdots  \lambda_i(M,g) \leq \hdots
\]
be the eigenvalues of the Laplacian.  For any constant $t > 0$,
$\lambda_i(M,tg)=\frac{\lambda_i(M,g)}{t}$. Let us consider  the
following functional which is invariant under dilations,
\[
\Lambda_i(M,g)=\lambda_i(M,g) Area(M,g)
\]
A result of Korevaar (\cite{18}) shows that there exists a constant
$C > 0$ such that for any $i > 0$ and surface $M$ of genus $\gamma$,
we have
\[
\Lambda_i(M,g) \leq C(\gamma +1)i
\]
For a given surface $M$ and a number $i$, one may ask what is the
supremum of $\Lambda_i(M,g)$ and which metric realizes it. It is a
very difficult question with just a few answers known, all but one
for $i=1$:
\begin{enumerate}
\item $\sup \Lambda_1(\mathbb S^2,g) = 8\pi$, and the maximum is the canonical metric on $\mathbb S^2$ (see \cite{19}).
\item $\sup \Lambda_1(\mathbb RP^2,g) = 12\pi$, and the maximum is the canonical metric on $\mathbb RP^2$ (see \cite{20}).
\item $\sup \Lambda_1(\mathbb T^2,g) = \frac{8\pi^2}{\sqrt{3}}$, and the maximum is the flat equilateral torus (see \cite{7}).
\item $\sup \Lambda_1(\mathbb K,g) = 12\pi E\left(\frac{2\sqrt{2}}{3}\right)$, and the maximum is the Lawson's bipolar surface
$\tilde \tau_{3,1}$ described below (see \cite{12}, \cite{13}). Here
$E(\alpha)$ is a complete elliptic integral of the second kind.
\item $\sup \Lambda_2(\mathbb S^2,g) = 16\pi$, and the maximum is a singular surface which is a union
of  two round spheres of area $2\pi$ glued together (see \cite{21}).
\end{enumerate}
The study of $\Lambda_i$-maximal metrics motivates the following
question. For every $i$, consider $\Lambda_i : g \to {\mathbb R}_+$ as
a functional on the space of Riemannian metrics on a surface. What
are the critical points of this functional? Functionals $\Lambda_i$
continuously depend on $g$ but in general are not differentiable.
However, for any analytic deformation $g_t$, $\lambda_i(M,g_t)$ is
left and right differentiable with respect to $t$ (see \cite{6}).
Following \cite{11} we propose the following definition for the
extremality of a metric:
\begin{definition}
If
\[
\frac{d}{dt} \lambda_i(M,g_t) \vert_{t=0^+} \leq 0 \leq \frac{d}{dt} \lambda_i(M,g_t) \vert_{t=0^-}
\]
for any analytic deformation $g_t$ preserving area and such that $g_0 = g$, the metric $g$ is said to be extremal for the
functional $\lambda_i$.
\end{definition}
Clearly, metrics (1)--(5) above are extremal since they are global
maxima.  An example of extremal metric for $\lambda_1(\mathbb
T^2,g)$ which is not a global maximum is given by the Clifford
torus. This metric together with (3) are the only critical metrics
for $\lambda_1$ on a torus (\cite{11}). Metrics (1), (2), (4) are
unique extremal metrics for $\lambda_1$ on a sphere, projective
plane and Klein bottle, respectively.

The following proposition proved in section \ref{section4} is an
extension to higher eigenvalues of what is shown for $\lambda_1$ in
\cite{11}.
\begin{prop}
\label{propo1} A metric $g$ is extremal for the functional
$\lambda_i$ if and only if there exists a finite family of
eigenfunctions $\phi_1,\dots ,\phi_m$ of $\lambda_i(M,g)$  with
$\sum_{k=1}^m d\phi_k \otimes d\phi_k = g$, and $j < i$ implies
$\lambda_j(M,g)<\lambda_i(M,g)$.
\end{prop}

Since extremal eigenvalues are always multiple \cite{11}, the last
condition of the Proposition means that extremality applies only to
the eigenvalue of the minimal rank.

It follows from Takahashi's theorem (see \cite{9}) that Proposition
\ref{propo1} is equivalent to the following statement:
$(\phi_1,\dots,\phi_m)$ is a minimal isometric immersion of $M$ into
a sphere ${\mathbb S}^{m-1}$ of radius
$\sqrt{\frac{2}{\lambda_i(g)}}$. In other words, any minimal surface
in a sphere carries an extremal metric for an eigenvalue of {\it
some} rank $i$.  Motivated by (4), in Theorem \ref{thm1} we compute
the ranks of extremal eigenvalues on Lawson's bipolar surfaces
$\tilde \tau_{r,k}$. As a byproduct, we find explicitly the
eigenfunctions providing the minimal immersion of a bipolar surface
into ${\mathbb S}^4$, see Theorem \ref{thm2}.

\smallskip

\noindent{\bf Remark.} It would be interesting to find also ranks of
the extremal eigenvalues on Lawson's surfaces themselves (they are
minimally immersed tori and Klein bottles in ${\mathbb S}^3$). Methods
developed in the present paper are not applicable straightforwardly
to this problem.

\subsection{Construction of bipolar surfaces}
The Lawson's surface $\tau_{r,k}$, with $ r > k > 0 $ and $(r,k) =
1$, is minimally immersed into $\mathbb S^3$ by
\[
I:\mathbb R^2 \rightarrow \mathbb R^4 \qquad I(u,v)=(\cos ru \cos v,\sin ru \cos v,\cos ku \sin v,\sin ku \sin v)
\]
It carries a metric  $(r^2\cos^2v+k^2\sin^2v) du^2 + dv^2$ and has a
group of symmetries that depends on $r$ and $k$.

If $r \equiv 0 \mod 2$ or $k \equiv 0 \mod 2$, $\tau_{r,k}$ is
invariant under the actions $(u,v)\rightarrow
(u+\pi,-v),(u,v)\rightarrow (u,v+2\pi)$ or $(u,v)\rightarrow
(u+\pi,\pi-v),(u,v)\rightarrow (u,v+2\pi)$ respectively. In those
cases, $\tau_{r,k}$ is a Klein bottle.

If $rk \equiv 1 \mod 2$, $\tau_{r,k}$ is a torus and its group of
motions is generated by $(u,v)\rightarrow (u+\pi,v+\pi)$ and
$(u,v)\rightarrow (u,v+2\pi)$.

Following \cite{14}, we realize the bipolar minimal surface $\tilde
\tau_{r,k}$ of $\tau_{r,k}$ as an exterior product of $I$ and
$I^{*}$, where $I^*$ is a unit vector normal to the torus (or Klein
bottle) $\tau_{r,k}$ and tangent to $\mathbb S^3$:
\begin{equation}
\label{itilde}
\tilde I = I \wedge I^{*} : \mathbb R^2 \rightarrow \mathbb S^5 \subset \mathbb R^6
\end{equation}
\[
I^{*}(u,v)=\frac{(k\sin ru \sin v, -k\cos ru \sin v, -r\sin ku \cos v,r\cos ku \cos v)}{\sqrt{r^2\cos^2v+k^2\sin^2v}}
\]
It will be shown in section \ref{equator} that $\tilde \tau_{r,k}$ actually lies in $\mathbb S^4$,
seen as an equator of $\mathbb S^5$.

\subsection{Spectral properties of bipolar surfaces}
\begin{theorem}
\label{thm1} For any $r,k \in \mathbb N$, with $0 < k < r$ and $(r,k) = 1$ we have the following,
\begin{enumerate}
\item If $rk \equiv 0 \mod 2$, $\tilde \tau_{r,k}$ is a torus and carries an extremal metric for the functional $\lambda_{4r - 2}$.
\item If $rk \equiv 1 \mod 4$, $\tilde \tau_{r,k}$ is a torus and carries an extremal metric for the functional $\lambda_{2r - 2}$.
\item If $rk \equiv 3 \mod 4$, $\tilde \tau_{r,k}$ is a Klein bottle and carries an extremal metric for the functional $\lambda_{r - 2}$.
\end{enumerate}
\end{theorem}
\noindent{\bf Remark.}  If $rk \equiv 3 \mod 4$, the double cover of
$\tilde \tau_{r,k}$, which is a torus, is an extremal metric for
$\lambda_{2r - 2}$.

\smallskip

For example, the Klein bottle $\tilde \tau_{5,3}$ and the double
cover of $\tilde \tau_{3,1}$ have extremal metrics for $\lambda_{3}$
and $\lambda_{4}$ respectively. Any other $\tilde \tau_{r,k}$ has an
extremal metric for $\lambda_i$ with $i \geq 5$ (for $\tilde
\tau_{2,1}, i = 6$).  It would be interesting to determine the type
of the critical metrics on bipolar surfaces. The only one that we
know is $\tilde \tau_{3,1}$ which is a global maximum.

In order to prove Theorem \ref{thm1} we need to analyze the
behaviour of the zeros of the  eigenfunctions providing the minimal
immersion of a bipolar surface.  Explicit expressions for the
eigenfunctions are given by
\begin{theorem}
\label{thm2} Set $n = r+k$,  $m = r-k$ if $rk \equiv 0 \mod 2$, and
$n = \frac{r+k}{2}$, $m = \frac{r-k}{2}$ if $rk \equiv 1 \mod 2$.
The minimal isometric immersion of a bipolar surface $\tilde
\tau_{r,k} \to {\mathbb S}^4 \subset {\mathbb R}^5$ is given by
\begin{equation}
\label{eigen} \left(\varphi_0(y), \cos(mx)\varphi_1(y),
\sin(mx)\varphi_1(y), \cos(nx)\varphi_2(y),
\sin(nx)\varphi_2(y)\right)\in {\mathbb S}^4,
\end{equation}
where $\varphi_0$, $\varphi_1$ and $\varphi_2$ are defined by the
following elliptic functions:
\begin{multline}
\label{wpfunctions}
\varphi_0(y)=\sqrt{\frac{n^2+m^2}{2n^2}}\left(1-\frac{n^2-m^2}{2\wp\left(y;a_{11},a_{12}\right)+b_1}\right)\\
\varphi_1(y)=\frac{1}{\sqrt{2}}\left(-1+\frac{n^2}{2\wp\left(y+\frac{1}{n}K\left(\frac{m}{n}\right);a_{21},a_{22}\right)+b_2}\right)\\
\varphi_2(y)=\sqrt{\frac{n^2-m^2}{2n^2}}\left(1+\frac{m^2}{2\wp\left(y;a_{31},a_{32}\right)+b_3}\right)\\
\end{multline}
and the matrices of constants are
\[
(a_{ij})=
\begin{pmatrix}
n^2m^2+\frac{(m^2+n^2)^2}{12} & & -\frac{n^2m^2(m^2+n^2)}{6} +\frac{(m^2+n^2)^3}{216}  \\
\\
m^2(m^2-n^2)+\frac{(2m^2-n^2)^2}{12} & & \frac{m^2(m^2-n^2)(2m^2-n^2)}{6}-\frac{(2m^2-n^2)^3}{216}\\
\\
n^2(n^2-m^2)+\frac{(2n^2-m^2)^2}{12} & & \frac{n^2(n^2-m^2)(2n^2-m^2)}{6}-\frac{(2n^2-m^2)^3}{216}\\
\end{pmatrix}
\]
\[
(b_i)=
\begin{pmatrix}
\frac{n^2-5m^2}{6}\\
\\
\frac{4m^2+n^2}{6}\\
\\
\frac{4n^2-5m^2}{6}\\
\end{pmatrix}
\]
\end{theorem}

\subsection{Plan of the paper} The paper is organized as follows.
In the next section, using a spectral-theoretic approach, we
construct ${\mathbb S}^1$-equivariant minimal tori and Klein bottles in
$\mathbb S^4$. In section 3 we show that these minimal surfaces and
the Lawson's bipolar surfaces are isometric, and prove Theorem
\ref{thm2}.  In section 4 we prove Theorem \ref{thm1} and find the
ranks of the extremal eigenvalues on Lawson's bipolar surfaces. In
the last section we compare our approach to constructing ${\mathbb
S}^1$-equivariant minimal tori in $\mathbb S^4$ with other methods
developed in \cite{15}, \cite{16}.


\section{Equivariant immersion of a torus in $\mathbb S^4$ by eigenfunctions}
\subsection{Basic facts}
\label{structure}
We consider an arbitrary torus $\mathbb T^2$ as a fundamental domain
in $\mathbb R^2$ for the group of motions generated by $(x,y)\rightarrow (x+2\pi,y)$,
$(x,y)\rightarrow (x,y+a)$, where $a > 0$ is a constant.  We suppose its metric conformal
and invariant under the actions $(x,y)\rightarrow (x+t,y)$, $(x,y)\rightarrow (x,-y)$,
$t \in \mathbb R$. It is given by
\[
f(y)(dx^2+dy^2)
\]
where $f(y) = f(-y) = f(y+a) > 0$.

The Laplacian on $\mathbb T^2$ is then
\[
\Delta = -\frac{1}{f(y)}\left(\frac{\partial^2}{\partial x^2}+\frac{\partial^2}{\partial y^2}\right)
\]
and the differential equation satisfied by the eigenfunction $\varphi$ of eigenvalue $\lambda$ is
\[
\Delta \varphi = \lambda \varphi
\]

\subsection{Immersion by the eigenfunctions}
The first step to prove theorem $\ref{thm2}$ is to construct a torus or Klein bottle $\rho_{n,m}$ admitting a minimal immersion by
its eigenfunctions (\ref{eigen}).
It is a simple generalization of what is done in \cite{12}.
We suppose the immersion is in $\mathbb S^4$, with radius $1$.
Takahashi's theorem (see \cite{9}) implies that these eigenfunctions are of eigenvalue $\lambda = 2$.
This gives the following conditions,
\begin{equation}
\varphi_0^2+\varphi_1^2+\varphi_2^2=1
\end{equation}
\begin{equation}
\label{metriccond}
(\varphi'_0)^2+(\varphi'_1)^2+(\varphi'_2)^2=m^2\varphi_1^2+n^2\varphi_2^2= f
\end{equation}

The $\varphi$-functions must then satisfy this system of second-order differential equations and all have the same period $a$,
\begin{equation}
\begin{cases}
\label{odesystem}
\varphi_0''=-2(m^2\varphi_1^2+n^2\varphi_2^2)\varphi_0\\
\varphi_1''=(m^2-2(m^2\varphi_1^2+n^2\varphi_2^2))\varphi_1\\
\varphi_2''=(n^2-2(m^2\varphi_1^2+n^2\varphi_2^2))\varphi_2\\
\end{cases}
\end{equation}

We consider only the case where $\varphi_1(y)$
is odd and both $\varphi_0(y)$ and $\varphi_2(y)$ are even functions.  This is because in this case we can find initial conditions
giving periodic solutions.
The following initial conditions
are required,
\begin{equation}
\label{condizero}
\begin{cases}
\varphi_0'(0)=0\\
\varphi_1(0)=0\\
\varphi_2'(0)=0\\
\end{cases}
\end{equation}

We can find two first integrals (see \cite{12} and \cite{13}) for the system (\ref{odesystem}),
\[
E_1 = (m^2 \varphi_1^2 + n^2 \varphi_2^2)^2 - (m^4 \varphi_1^2 + n^4 \varphi_2^2) + m^2 (\varphi_1^{'})^2 + n^2  (\varphi_2^{'})^2
\]
\[
E_2 = n^2 (n^2-m^2) \varphi_2^2 (\varphi_2^2 -1)
+ m^2(n^2-m^2) \varphi_2^2 \varphi_1^2
+ m^2 \varphi_2^2 (\varphi_1')^2
\]
\[
\qquad \qquad \qquad \qquad \qquad \qquad \qquad \qquad -2 m^2 \varphi_1 \varphi_2 \varphi_1' \varphi_2'
+ (\varphi_2')^2 ((n^2-m^2) + m^2 \varphi_1^2)
\]
Since (\ref{metriccond}) and (\ref{condizero}) give $\varphi_1'(0)^2=n^2\varphi_2^2(0)$, the expression for $E_1$ reduces to
\[
E_1(0) = n^4 \varphi_2^4(0) - n^2(n^2-m^2) \varphi_2^2(0)
\]
It admits a minimum for $\varphi_2^2(0)=\frac{n^2-m^2}{2n^2}$, which corresponds to a periodic solution. Indeed if we set
\begin{equation}
\begin{cases}
\varphi_0(0)=\sqrt{\frac{n^2+m^2}{2n^2}}\\
\varphi_1'(0)=\sqrt{\frac{n^2-m^2}{2}}\\
\varphi_2(0)=\sqrt{\frac{n^2-m^2}{2n^2}}\\
\end{cases}
\end{equation}
the solution lies on the intersection of the sphere and the cylinder of equation
\[
2\varphi_1^2+\frac{2n^2}{n^2+m^2}\varphi_0^2 = 1
\]
By setting
\[
\varphi_0(y)=\sqrt{\frac{n^2+m^2}{2n^2}}\cos\theta(y) \qquad \varphi_1(y)=\frac{1}{\sqrt{2}}\sin\theta(y)
\]
we have the following differential equation and solution for $\theta(y)$,
\[
(\theta')^2 = n^2-m^2\cos^2\theta
\]
\[
y = \frac{1}{n} \int_{0}^{\theta} \frac{d\theta}{\sqrt{1-\left(\frac{m}{n}\right)^2\cos^2\theta}}
\]
\newline
Using this solution we find
\begin{equation}
\label{perioda}
a = \frac{1}{n} \int_{0}^{2\pi} \frac{d\theta}{\sqrt{1-\left(\frac{m}{n}\right)^2\cos^2\theta}} = \frac{4}{n} K\left(\frac{m}{n}\right)
\end{equation}
where $K\left(\frac{m}{n}\right)$ is a complete elliptic integral of the first kind.

In particular, we have $\varphi_2(y) \neq 0$ for any $y$ and the following metric
\begin{equation}
\label{normmetric}
2 f(y)(dx^2+dy^2) = \left[ (m^2+n^2) - 2m^2\cos^2\theta(y)\right](dx^2+dy^2)
\end{equation}

After some substitutions, the ODEs of (\ref{odesystem}) can be integrated to become the following uncoupled system of first-order,
\begin{equation}
\begin{cases}
(\varphi_0')^2=\frac{2n^2m^2}{n^2+m^2}\varphi_0^4-(m^2+n^2)\varphi_0^2+\frac{m^2+n^2}{2}\\
(\varphi_1')^2=-2m^2\varphi_1^4+(2m^2-n^2)\varphi_1^2+\frac{n^2-m^2}{2}\\
(\varphi_2')^2=-2n^2\varphi_2^4+(2n^2-m^2)\varphi_2^2+\frac{m^2-n^2}{2}\\
\end{cases}
\end{equation}
Each of these equations has a solution in terms of the Weierstrass $\wp$ function, given by (\ref{wpfunctions}).
If $n$ is odd or $m$ is even, $\rho_{n,m}$ is the torus with stucture presented in section \ref{structure}, where the period $a$ and the metric $f(y)$ are given by (\ref{perioda}) and (\ref{normmetric}) respectively.

In the case $n$ is even and $m$ is odd, we may consider the domain of the eigenfunctions to be a Klein
bottle $\mathbb K$ instead of a torus.  Indeed, we can then enlarge the group of motions by adding the generator
$(x,y)\rightarrow (x+\pi,-y)$ to the ones described in \ref{structure}. The eigenfunctions are still well-defined because of the
parity of the $\varphi$-functions.
When $n$ is even and $m$ is odd, $\rho_{n,m}$ is the Klein bottle with metric $f(y)$ in (\ref{normmetric}) and period $a$ in (\ref{perioda}).

\smallskip
\noindent{\bf Remark.}  As in \cite{13}, the first integral $E_1$ can be considered as an Hamiltonian $H(q_1,q_2,p_1,p_2)$ by setting,
\[
q_1 = \varphi_1 \qquad p_1 = 2m^2 \varphi_1^{'}
\]
\[
q_2 = \varphi_2 \qquad p_2 = 2n^2 \varphi_2^{'}
\]
Since $E_2$ is a second independent first integral, the equations for $\varphi_1$ and $\varphi_2$ in (\ref{odesystem})
form an integrable system. We can also find $E_2$ from $E_1$ by writing $(\varphi_0\varphi_0^{'})^2$ in the two following ways
\[
(\varphi_1\varphi'_1 + \varphi_2\varphi_2^{'})^2 = (1-\varphi_1^2-\varphi_2^2)(m^2\varphi_1^2+n^2\varphi_2^2 - (\varphi'_1)^2 - (\varphi'_2)^2)
\]


\section{Bipolar surfaces of the Lawson's tori or Klein bottles}
\subsection{Properties of $\tilde \tau_{r,k}$}
\label{equator}
The explicit parametric representation of (\ref{itilde}) is given by,
\[
\tilde I(u,v)= \frac{1}{\sqrt{r^2\cos^2v+k^2\sin^2v}}
\begin{pmatrix}
-k \sin v \cos v\\
r \sin v \cos v\\
-r \cos^2v \sin ku \cos ru - k \sin^2v \sin ru \cos ku\\
r \cos^2v \sin ru \cos ku + k \sin^2v \sin ku \cos ru\\
-r \cos^2v \sin ku \sin ru + k \sin^2v \cos ru \cos ku\\
r \cos^2v \cos ku \cos ru - k \sin^2v \sin ru \sin ku\\
\end{pmatrix}
\]
\newline
If we apply the following orthogonal transformation to this vector we get a cleaner result,
\[
A = \frac{1}{\sqrt{2}}
\begin{pmatrix}
1 & 1 & & &0 &\\
-1 & 1 & & & &\\
 &  & 1 & 1 & &\\
 &  & -1 & 1 & &\\
 &  & & & 1 & 1\\
 &  0& & & -1 & 1\\
\end{pmatrix}
\]
\newline
\begin{equation}
\label{taucie}
A \circ \tilde I(u,v) = \frac{1}{\sqrt{8}\sqrt{r^2\cos^2v+k^2\sin^2v} }
\begin{pmatrix}
(r - k) \sin 2v\\
(r + k) \sin 2v\\
[(r-k) + (r+k) \cos 2v] \sin (r-k)u\\
[(r+k) + (r-k) \cos 2v] \sin (r+k)u\\
[(r+k) + (r-k) \cos 2v] \cos (r+k)u\\
[(r-k) + (r+k) \cos 2v] \cos (r-k)u\\
\end{pmatrix}
\end{equation}
Since the image of $A \circ \tilde I$ is always orthogonal to the vector $(r+k,k-r,0,0,0,0)$, we may consider it to lie
in $\mathbb S^4 \subset \mathbb R^6$ instead of $\mathbb S^5$.

The bipolar surface immersed by (\ref{taucie}) has the metric
\newline
\begin{equation}
\frac{(r^2 - (r^2-k^2)\sin^2v)^2 + r^2k^2}{r^2 - (r^2-k^2)\sin^2v} \left(du^2 + \frac{dv^2}{r^2 - (r^2-k^2)\sin^2v} \right)
\end{equation}

If $rk \equiv 0 \mod 2$, the group of motions of the bipolar surface
is generated by $(u,v)\rightarrow (u,v + \pi)$, $(u,v)\rightarrow (u + 2\pi,v)$ in $\mathbb R^2$, which corresponds to a torus.  Thus if $\tau_{r,k}$
is a Klein bottle, $\tilde \tau_{r,k}$ is a torus. Furthermore, we can set $n = r+k$ and $m = r-k$ to write the previous metric as
\begin{equation}
\label{gmetric}
g_{0} = \frac{\left[\left(n+m\right)^2-4mn\sin^2v\right]^2+\left(n^2-m^2\right)^2}{\left(n+m\right)^2-4mn\sin^2v}\left(\frac{du^2}{4}+
\frac{dv^2}{\left(n+m\right)^2-4mn\sin^2v}\right)
\end{equation}

\begin{lemma}
\label{previouslemma}
If $rk \equiv 1 \mod 4$, $\tilde \tau_{r,k}$is a torus. If $rk \equiv 3 \mod 4$, $\tilde \tau_{r,k}$ is a Klein
bottle.
\end{lemma}
\begin{proof}
If $rk \equiv 1 \mod 2$, we can set $n = \frac{r+k}{2}$ and $m = \frac{r-k}{2}$.
The previous parametric representation and metric of the bipolar surface then become,
\begin{equation}
\label{parambip}
A \circ \tilde I(u,v) = \frac{1}{\sqrt{2}\sqrt{\left(n+m\right)^2-4mn\sin^2v} }
\begin{pmatrix}
m \sin 2v\\
n \sin 2v\\
[m + n \cos 2v] \sin 2mu\\
[n + m \cos 2v] \sin 2nu\\
[n + m \cos 2v] \cos 2nu\\
[m + n \cos 2v] \cos 2mu\\
\end{pmatrix}
\end{equation}
\begin{equation}
\label{metric2}
g_{0}' = \frac{\left[\left(n+m\right)^2-4mn\sin^2v\right]^2+\left(n^2-m^2\right)^2}{\left(n+m\right)^2-4mn\sin^2v}\left(du^2+
\frac{dv^2}{\left(n+m\right)^2-4mn\sin^2v}\right)
\end{equation}
If $rk \equiv 1 \mod 4$, then $n$ is odd and $m$ is even.
If we study the group of invariant transformations for (\ref{parambip}),
we find it to be generated by $(u,v)\rightarrow (u+\pi,v)$, $(u,v)\rightarrow (u,v+\pi)$ in $\mathbb R^2$.
This corresponds to a torus.

Now if $rk \equiv 3 \mod 4$, $n$ is even and $m$ is odd so that we may add the generator
$(u,v)\rightarrow H_1^{-1} \circ H_2 \circ H_1(u,v)$ to the group and turn the torus in a Klein bottle. The following
transformations are used,
\[
H_1(u,v) = \left(u,\frac{1}{n+m}\int_{0}^{v}\frac{dv}{\sqrt{1-\left(\frac{2\sqrt{mn}}{n+m}\right)^2\sin^2v}}\right) = (u,z)
\]
\[
H_2(u,z)=\left(u+\frac{\pi}{2}, \frac{a}{4} - z\right)
\]
This can be verified directly by using some identities, for $\alpha = \frac{m}{n}$ and $\alpha' = \sqrt{\frac{n^2-m^2}{n^2}}$.
First, we use the following expressions for $\sin v$ and $\cos v$ in (\ref{parambip}) and (\ref{metric2}) to apply $H_1$,
\[
\sin v = sn\left[(1+\alpha) nz, \frac{2\sqrt{\alpha}}{1+\alpha}\right] = \left(1+\alpha\right)\frac{sn\left(nz,\alpha\right)}{1+\alpha sn^2\left(nz,\alpha\right)}
\]
\[
\cos v = cn\left[(1+\alpha) nz, \frac{2\sqrt{\alpha}}{1+\alpha}\right] = \frac{cn\left(nz,\alpha\right)dn\left(nz,\alpha\right)}{1+\alpha sn^2\left(nz,\alpha\right)}
\]
\[
\sqrt{\left(1+\alpha\right)^2-4\alpha \sin^2v} = \left(1 + \alpha\right) \frac{1 - \alpha sn^2\left(nz,\alpha\right)}{1+ \alpha sn^2\left(nz,\alpha\right)}
\]
\newline
After applying $H_2$, we use the identities below to get back from $z$ to $v$ dependence with $H_1^{-1}$, and show
that (\ref{parambip}) and (\ref{metric2}) are invariant under $(u,v)\rightarrow H_1^{-1} \circ H_2 \circ H_1(u,v)$.
\[
\frac{1}{m+n} K\left(\frac{2\sqrt{\alpha}}{1+\alpha}\right) = \frac{1}{n} K\left(\alpha\right) = \frac{a}{4}
\]
\[
sn(-w,\alpha) = - sn(w,\alpha) \qquad cn(-w,\alpha) = cn(w,\alpha) \qquad dn(-w,\alpha) = dn(w,\alpha)
\]
\[
sn(w+K(\alpha),\alpha) = \frac{cn(w,\alpha)}{dn(w,\alpha)} \qquad cn(w+K(\alpha),\alpha) = -\alpha' \frac{sn(w,\alpha)}{dn(w,\alpha)}
\]
\[
dn(w+K(\alpha),\alpha) = \alpha' \frac{1}{dn(w,\alpha)}
\]
\[
sn^2(w,\alpha)+cn^2(w,\alpha) = 1 \qquad dn^2(w,\alpha)+\alpha ^2 sn^2(w,\alpha) = 1
\]
\end{proof}

\subsection{Relation between $\tilde \tau_{r,k}$ and $\rho_{n,m}$}
\begin{theorem}
\label{thm3}
If $rk \equiv 0 \mod 2$, $\tilde \tau_{r,k}$ is isometric to $\rho_{r+k,r-k}$.
Otherwise, $\tilde \tau_{r,k}$ is isometric to $\rho_{\frac{r+k}{2},\frac{r-k}{2}}$.
\end{theorem}
\begin{proof}
We need to define two other transformations
\[
H_3(u,z) = \left(u,2z + \frac{1}{n}K\left(\frac{m}{n}\right)\right) = (x,y)
\]
\[
H_{3}'(u,z) = \left(2u,2z + \frac{1}{n}K\left(\frac{m}{n}\right)\right) = (x,y)
\]
We recall that the metric on $\rho_{n,m}$ (see (\ref{normmetric}) ) is given by
\[
g(n,m)=\left[\frac{m^2+n^2}{2} - m^2\cos^2\theta(y)\right](dx^2+dy^2)
\]
When $rk \equiv 0 \mod 2$, we set $n = r+k$ and $m = r-k$. The metric $g_{0}$ on $\tilde \tau_{r,k}$ is given by (\ref{gmetric}) and we have
\[
H_1^{*} \circ H_3^{*} g(n,m) = g_{0}
\]
When $rk \equiv 1 \mod 2$, we set $n = \frac{r+k}{2}$ and $m = \frac{r-k}{2}$. The metric $g_{0}'$ on $\tilde \tau_{r,k}$ is given by (\ref{metric2}) and we have
\[
H_1^{*} \circ H_{3}^{'*} g(n,m) = g_{0}'
\]
This can be verified with the additional identities,
\[
\cos\theta (y) = sn\left[K\left(\alpha\right) - ny, \alpha\right]
\]
\[
\cos\theta (z) = -sn\left(2nz,\alpha\right) = -\frac{2sn\left(nz,\alpha\right) cn\left(nz,\alpha\right) dn\left(nz,\alpha\right)}{1-\alpha^2 sn^4\left(nz,\alpha\right)}
\]
\[
\sin\theta (z) = cn\left(2nz,\alpha\right) = \frac{1 - 2sn^2\left(nz,\alpha\right) + \alpha^2 sn^4\left(nz,\alpha\right)}{1-\alpha^2 sn^4\left(nz,\alpha\right)}
\]
Using lemma \ref{previouslemma} and previous considerations, it is easy to verify that the groups of
motions in $\mathbb R^2$ correspond to each other, whether the surface is a torus or a Klein bottle.
We may then conclude that our surfaces are isometric.
\end{proof}
Theorem \ref{thm2} follows from the results of Section 2 and Theorem \ref{thm3}.
Note that each coordinate function in (\ref{parambip}) is one of the functions in (\ref{eigen}) multiplied by a scalar.


\section{The rank of the extremal eigenvalue of bipolar surfaces}
\subsection{The spectrum of $\rho_{n,m}$}

Let the sequence
\[
0 = \lambda_0(g) < \lambda_1(g) \leq \lambda_2(g) \leq \lambda_3(g) \leq \hdots \leq \lambda_k(g) \leq \hdots
\]
be the eigenvalues of the Laplacian on a torus (or Klein bottle) represented by $g$, counting multiplicities.
Theorem \ref{thm3} implies that $\tilde \tau_{r,k}$ and $\rho_{n,m}$ have the same eigenvalues for the right choice of $n$ and $m$.
We will use the $(x,y)$ coordinates of the corresponding $\rho_{n,m}$ to study the spectrum of $\tilde \tau_{r,k}$.

Since the Laplacian $\Delta$ and the operator $\frac{\partial}{\partial x}$ on $\rho_{n,m}$ commute in the (x,y) coordinates,
the eigenspaces of $\rho_{n,m}$ admit functions of the type $\sin(px)\phi(y)$ and $\cos(px)\phi(y)$ as basis.
This gives the following eigenvalue problem for $\phi$ with periodic boundary conditions $\phi(y) = \phi(y+a)$:
\begin{equation}
\label{odespec}
\phi^{''} + \left[\lambda f(y) - p^2\right]\phi = 0
\end{equation}
We may then study the spectrum of the above equation for each $0 \leq p \in \mathbb N$ to get the complete spectrum of $\rho_{n,m}$.
Note that $f(y)$ is just the metric of $\rho_{n,m}$, i.e.
\[
f(y) = \frac{m^2+n^2}{2} - m^2\cos^{2}\theta (y)
\]
with $f(y) = f(-y) = f(y + \frac{a}{2}) > 0$.

Haupt's theorem states that, for each $p$, there exists a sequence of eigenvalues for (\ref{odespec}), counting multiplicities :
\[
\gamma_0(p) < \gamma_1(p) \leq \gamma_2(p) < \hdots < \gamma_{2n-1}(p) \leq \gamma_{2n}(p) < \hdots
\]

The respective eigenfunctions $\Phi_0$ have no zeros and both $\Phi_{2n-1}$ and $\Phi_{2n}$ have 2n zeros on $[0,a)$.
We know that $\gamma_0(0)=0$ and that $\gamma_0(p)>0$ if $p > 0$.  Since $\varphi_2(y)$ has no zeros, $\gamma_0(n)=2$.

\subsection{Coexistence problem} The first step to prove theorem \ref{thm1} is to study the multiplicities in (\ref{odespec}) for a fixed $p$.
\begin{lemma}
\label{lemmacoex}
If $ 0 < \gamma_i(p) < 3 $, then $mult(\gamma_i) = 1$ as an eigenvalue of (\ref{odespec}) with fixed parameter $p$.
\end{lemma}
\begin{proof}
If we express the differential equation in the variable $\theta$,  we get the following Ince's equation with boundary condition
$\phi(\theta) = \phi(\theta + 2\pi)$ :
\[
(1 + c_{1} \cos 2\theta ) \phi^{''}(\theta) + c_{2} \sin(2\theta) \phi^{'}(\theta) + ( c_{3}+ c_{4}\cos 2\theta) \phi(\theta) = 0
\]
\begin{equation}
\begin{cases}
c_{1} = -\frac{\alpha^2}{2-\alpha^2}\\
c_{2} = -2c_{1}\\
c_{3} = \frac{\lambda - 2\left(\frac{p}{n}\right)^2 }{2-\alpha^2}\\
c_{4} = \lambda c_{1}\\
\end{cases}
\end{equation}
with $|c_{1}| < 1$ ($\alpha = \frac{m}{n}$).

A theorem about this equation (see \cite{3} theorem 7.1) states that if there are 2 linearly independent $\pi$-periodic or $2\pi$-periodic solutions associated to the eigenvalue $\lambda$, then the
polynomial $Q(\mu) = \mu^2 + 2\mu -\lambda$ has an integral root. This means that there exists an integer $c$ such that
$\lambda = c^2 -1$.  This is impossible if $ 0 < \gamma_i(p) < 3 $, so $mult(\gamma_i) = 1$.
\end{proof}

In particular, for $p = 0,m,n$, the $\varphi$-functions in
(\ref{wpfunctions}) are the only eigenfunctions of eigenvalue $2$.

\subsection{Monotonicity of the spectrum}  Set $R = p^2$, $b = \frac{a}{2}$ and let $z_1(y;R,\lambda)$, $z_2(y;R,\lambda)$ be the solutions to
(\ref{odespec}) (without the boundary conditions) with the following initial conditions:

\begin{equation}
\begin{cases}
z_{1}(0)=1\\
z_{1}^{'}(0)=0\\
z_{2}(0)=0\\
z_{2}^{'}(0)=1\\
\end{cases}
\end{equation}
We define $\Psi(R,\lambda)$ = $z_1(b;R,\lambda)$ + $z_{2}^{'}(b;R,\lambda)$ which we assume to be smooth.
The oscillation theorem states that $\lambda$ is
an eigenvalue of (\ref{odespec}) if and only if $\Psi(R,\lambda)^2 = 4$.

Its proof uses the following identities to find solutions $g(y) = k_1 z_1(y) + k_2 z_2(y)$ such that $g(y+b)=C g(y)$ for an undetermined
constant C.
\begin{equation}
\label{periodswitch}
\begin{cases}
z_{1}(y+b) = z_{1}(b) z_{1}(y) + z_{1}^{'}(b) z_{2}(y)\\
z_{2}(y+b) = z_{2}(b) z_{1}(y) + z_{2}^{'}(b) z_{2}(y)\\
\end{cases}
\end{equation}
It follows that $C$ must be a root of the characteristic equation $P(x) = x^2 - \Psi x + 1$, (note that $z_{1}(y)z_{2}^{'}(y) - z_{2}(y)z_{1}^{'}(y) = 1$).
If $\Psi^2 \neq 4$, there are two independent solutions $g_1(y),g_2(y)$ associated to distinct constants $C_{1},C_{2}$.  Any linear
combination of them is not a solution of period $2b$ since this would imply that $C_{1}^{2}=C_{2}^{2}=1$ and contradict the fact that
$C_{1}C_{2}=1$ and $C_{1} \neq C_{2}$.
But if $\Psi^2 = 4$, it follows that $C^2=1$ and then there exists a periodic solution $g(y+2b)=C^2g(y)=g(y)$.

\begin{lemma}
The functions $\gamma_i(p)$ are strictly increasing in $p$ if $p \geq 0$ and $0 < \gamma_i(p) < 3$.
\end{lemma}
\begin{proof}
We want to study locally the $\gamma$-functions defined for  $0 < p \in \mathbb R$, such that $\Psi(p^2,\gamma_i(p))^2 = 4$.
To use the implicit function theorem, we must first verify that $\frac{\partial \Psi}{\partial \lambda} \neq 0$ at any point $(p^2,\gamma_i(p))$.
According to the oscillation theorem, $\frac{\partial \Psi}{\partial \lambda}(R,\lambda) = 0$ and $\Psi(R,\lambda)^2 = 4$
implies that $\lambda$ has multiplicity 2. We've shown in lemma \ref{lemmacoex} that this is impossible if $0 < \lambda < 3$.
We conclude that the $\gamma$-functions can be expressed locally as a function of $p$ in the domain concerned.

The derivatives of $\Psi(R,\lambda)$ are
\[
\frac{\partial \Psi}{\partial \lambda} = [z_{1}(b)-z_{2}^{'}(b)]\int_{0}^{b} f(y) z_1 z_2 dy + z_{1}^{'}(b) \int_{0}^{b} f(y) z_2^2 dy
 - z_{2}(b) \int_{0}^{b} f(y) z_1^2 dy
\]
\[
\frac{\partial \Psi}{\partial R} = [z_{2}^{'}(b)-z_{1}(b)]\int_{0}^{b} z_1 z_2 dy - z_{1}^{'}(b) \int_{0}^{b} z_2^2 dy
 + z_{2}(b) \int_{0}^{b} z_1^2 dy
\]
But $f(y) = f(-y)$ gives $z_1(y) = z_1(-y)$, $z_2(y) = -z_2(-y)$ and
\begin{equation}
\begin{cases}
z_{1}(b) = 2z_{1}(b/2)z_{2}^{'}(b/2) -1 = 1 + 2z_{2}(b/2)z_{1}^{'}(b/2)\\
z_{1}^{'}(b) = 2z_{1}(b/2)z_{1}^{'}(b/2)\\
z_{2}(b) = 2z_{2}(b/2)z_{2}^{'}(b/2)\\
z_{2}^{'}(b) = z_{1}(b)\\
\end{cases}
\end{equation}
These identities can be verified using (\ref{periodswitch}) and their derivatives, by putting $y = -b/2$ and solving the linear system obtained.
The derivatives of $\Psi$ can then be written
\[
\frac{\partial \Psi}{\partial \lambda} = z_{1}^{'}(b) \int_{0}^{b} f(y) z_2^2 dy - z_{2}(b) \int_{0}^{b} f(y) z_1^2 dy
\]
\[
\frac{\partial \Psi}{\partial R} = - z_{1}^{'}(b) \int_{0}^{b} z_2^2 dy  + z_{2}(b) \int_{0}^{b} z_1^2 dy
\]

The Wronskian of the system (\ref{odespec}) is a constant $z_{1}(y)z_{2}^{'}(y) - z_{2}(y)z_{1}^{'}(y) = 1$.
Since $z_{2}^{'}(b) = z_{1}(b)$ and $\Psi(R,\lambda)^2=(z_{2}^{'}(b) + z_{1}(b))^2=4$, we have $z_{1}(b)z_{2}^{'}(b) = 1$ and $z_{2}(b)z_{1}^{'}(b) = 0$.
There is exactly one of the terms $z_{2}(b),z_{1}^{'}(b)$ that is zero (if it was both, $\frac{\partial \Psi}{\partial \lambda}$ would vanish).
In any case we get from the implicit function theorem
\[
\frac{d \gamma_i(p)}{d(p^2)} =
\frac{d \lambda}{dR} = - \frac{\frac{\partial \Psi}{\partial R}}{\frac{\partial \Psi}{\partial \lambda}} > 0
\]
if $ p > 0$ and conclude that the functions $\gamma_i(p)$ are strictly increasing as functions of $p$ if $p \geq 0$ and
$0 < \gamma_i(p) < 3$.
\end{proof}

\subsection{Multiplicity and parity of the eigenfunctions}
\label{section4}
Both $\varphi_0(y)$ and $\varphi_1(y)$ have two zeros
on their period $a$.  This restricts their rank to $\gamma_1(p)$ or $\gamma_2(p)$
in the sequences of eigenfunctions. Since the $\gamma$-functions are strictly increasing, it must be
$\gamma_2(0)=2$ (for $\varphi_0(y)$) and $\gamma_1(m)=2$ (for $\varphi_1(y)$). The multiplicity of $\lambda = 2$ is then $5$.

At each point $(p^2,\gamma_i(p))$, either $z_{1}^{'}(b) = 0$ or $z_{2}(b) = 0$. We can study the parity of the eigenfunctions using (\ref{periodswitch}).
If $z_{1}^{'}(b) = 0$, the eigenfunction is $z_{1}(y)$ and even.
If $z_{2}(b) = 0$, the eigenfunction is $z_{2}(y)$ and odd.
The eigenfunctions associated to $\gamma_0(p)$ have no zeros, so they must all be even.

Since $\varphi_1(y)$ is odd, we know $z_{2}(b) = 0$ at $p = m$. We can extend this
equality to all $0 \leq p \leq m$. Indeed,
\[
A = \left\{ p \in [0,m] | z_{1}^{'}(b;p^2,\gamma_1(p)) \neq 0 \right\} = \left\{ p \in [0,m] | z_{2}(b;p^2,\gamma_1(p)) = 0 \right\}
\]
because $0 < \gamma_1(p) < 3$ for $p \in [0,m]$ and $z_{2}(b)$ and $z_{1}^{'}(b)$ cannot vanish at the same time.
The set $A$ is then non-empty, open and closed in $[0,m]$, so it has to be the whole interval.
The eigenfunctions associated to $\gamma_1(p)$ must then be odd for $0 \leq p \leq m$.

\begin{lemma}
\label{lastlemma}
On the surface $\rho_{n,m}$, $mult(2) = 5$.
If $\rho_{n,m}$ is a torus, $\lambda_{i} = 2$ if and only if $2(n+m-1) \leq i \leq 2(n+m+1)$.
If $\rho_{n,m}$ is a Klein bottle, $\lambda_{i} = 2$ if and only if $(n+m-2) \leq i \leq (n+m+2)$.
\end{lemma}
\begin{proof}
We've just shown $mult(2) = 5$. To find the rank of $\lambda = 2$, we count $2(n-1)$ non-zero eigenvalues
smaller than 2 for each $\gamma_0(p)$, $ 0 \leq p \leq n-1$.
The factor $2$ is the multiplicity from $\sin(px)\phi(y)$ and $\cos(px)\phi(y)$. For $\gamma_1(p)$, there are $2(m-1) + 1$.  This gives
a total of $2(n+m) - 3$ non-zero eigenvalues smaller than $2$ when $\rho_{n,m}$ is a torus.

If $\rho_{n,m}$ is a Klein bottle, $m$ is odd and $n$ is even.  We must then reject the eigenvalues $\gamma_0(p)$ when $p$ is odd
(because $\phi(y)$ is even for $\gamma_0$) and $\gamma_1(p)$ when $p$ is even (because $\phi(y)$ is odd for $\gamma_1$). This gives $n+m-3$
non-zero eigenvalues smaller than $2$.
\end{proof}

\begin{proof}[Proof of proposition $\ref{propo1}$]
If $\lambda$ is an eigenvalue of $(M,g)$ with multiplicity $r$,
lemma 3.15 in \cite{6} states that for any analytic deformation $g_t$ such that $g_0 = g$, there exists $r$ scalars $\Pi_j$
and $r$ functions $\Phi_j$ depending on $t$ such that:
\begin{enumerate}
\item $\Delta \Phi_j = \Pi_j \Phi_j$ for any $j$ and $t$,
\item $\Pi_j(0)= \lambda$ for any $j$,
\item {$\Phi_j$} is orthonormal for any $t$.
\end{enumerate}
If $i$ is the smallest integer such that $\lambda_i(M,g) = \lambda$, then $g$ will be an extremal metric of $\lambda_i$ if and only if
for any analytic deformation $g_t$ there are both nonpositive and nonnegative values for $\Pi_{j}^{'}(0)$.
Indeed, we'll have
\[
\frac{d}{dt} \lambda_i(M,g_t) \vert_{t=0^+} = \min_{j} \Pi_{j}^{'}(0) \leq 0 \leq \max_{j} {\Pi_{j}^{'}(0)}  = \frac{d}{dt} \lambda_i(M,g_t) \vert_{t=0^-}
\]
An addendum to \cite{11} shows that if there exists a finite family \{$\phi_k$\} of eigenfunctions of $\lambda$ on $(M,g)$
with $\sum d\phi_k \otimes d\phi_k = g$, then for any deformation $g_t$
\[
\min_{j} \Pi_{j}^{'}(0) \leq 0 \leq \max_{j} {\Pi_{j}^{'}(0)}
\]
This implies the extremality of $g$ for the functional $\lambda_i$, if $j < i$ implies $\lambda_j(M,g)<\lambda_i(M,g) = \lambda$.
\end{proof}

\begin{proof}[Proof of theorem $\ref{thm1}$]
Since $\tilde \tau_{r,k}$ admits an isometric immersion in $\mathbb S^4$ (theorems \ref{thm2} and \ref{thm3}) by eigenfunctions of $\lambda = 2$,
the last proposition implies that $\tilde \tau_{r,k}$ is an extremal metric for some functional $\lambda_i$.
The number $i$ will be the smallest integer such that $\lambda_i = 2$ on $\tilde \tau_{r,k}$ and we find it using lemma \ref{lastlemma}.
\begin{enumerate}
\item If $rk \equiv 0 \mod 2$, $i = 2((r+k)+(r-k)-1) = 4r - 2$.
\item If $rk \equiv 1 \mod 4$, $i = 2(\frac{r+k}{2}+\frac{r-k}{2}-1) = 2r - 2$.
\item If $rk \equiv 3 \mod 4$, $i = \frac{r+k}{2}+\frac{r-k}{2}-2 = r - 2$.
\end{enumerate}
\end{proof}

\noindent{\bf Remark.}
We may calculate the value of $\Lambda_i(\tilde \tau_{r,k})$, where $i$ is the rank of the extremal eigenvalue given in the
preceding theorem.
We first calculate the area of a torus $\rho_{n,m}$ as we did for the period $a$ in (\ref{perioda}), using the $\theta$
change of variable.
\[
Area(\rho_{n,m}) = 2\pi\int_{0}^{a} f(y) dy = \pi\int_{0}^{a} \left((m^2+n^2)-2m^2\cos^2\theta(y)\right) dy =
\]
\begin{equation}
\label{lambdaarea}
4\pi n \left[\left(\frac{m^2}{n^2}-1\right)K\left(\frac{m}{n}\right)+2E\left(\frac{m}{n}\right)\right]=
\end{equation}
\[
4\pi (n+m) E\left(\frac{2\sqrt{mn}}{m+n}\right)
\]
where $E(\frac{m}{n})$ is a complete elliptic integral of the second kind.
Note that if $\rho_{n,m}$ is a Klein bottle, the area must be divided by two to get $Area(\rho_{n,m})$.
Using the correspondence between the surfaces $\rho_{n,m}$ and $\tilde \tau_{r,k}$ and the fact that $\lambda=2$
in the case we are interested in, we obtain
\begin{enumerate}
\item If $rk \equiv 0 \mod 2$, $\Lambda_{4r - 2}(\tilde \tau_{r,k})= 16\pi r E\left(\frac{\sqrt{r^2-k^2}}{r}\right)$.
\item If $rk \equiv 1 \mod 4$, $\Lambda_{2r - 2}(\tilde \tau_{r,k})= 8\pi r E\left(\frac{\sqrt{r^2-k^2}}{r}\right)$.
\item If $rk \equiv 3 \mod 4$, $\Lambda_{r - 2}(\tilde \tau_{r,k})= 4\pi r E\left(\frac{\sqrt{r^2-k^2}}{r}\right)$.
\end{enumerate}


\section{$\mathbb S^1$-equivariant maps into spheres and minimal surfaces}
\subsection{Minimal tori in $\mathbb S^4$}
The aim of this section is to compare our approach to study minimal immerions in $\mathbb S^4$ with those previously used in
\cite{15} and \cite{16}.
The minimal immersion given by (\ref{eigen}) is rotationally symmetric in $\mathbb S^4$ for the $\mathbb S^1$-action represented
by the matrices
\begin{equation}
R(\theta) =
\begin{pmatrix}
1 & & & 0&\\
 & \cos(m\theta)& -\sin(m\theta) & & \\
 & \sin(m\theta)&  \cos(m\theta)& & \\
 &  &  & \cos(n\theta) & -\sin(n\theta)\\
 &  0& & \sin(n\theta)& \cos(n\theta) \\
\end{pmatrix}
\end{equation}

Immersions in $\mathbb S^4$ equivariant under the rotations $R(\theta)$ are studied in \cite{15}, using the integrable system approach.
Minimal surfaces then correspond to geodesics in the orbit space of $\mathbb S^4$ under $R(\theta)$, which is parameterized
by the variables $(\rho,\alpha,\psi)$. This space has the following metric and identification in $\mathbb S^4$,
\[
G = \cos^2\rho (m^2 \cos^2\alpha + n^2 \sin^2\alpha) (d\rho^2 + \cos^2\rho d\alpha^2) + m^2 \sin^2\alpha \cos^2\alpha \cos^4\rho d\psi^2
\]
\[
\begin{pmatrix}
\sin\rho \\
\cos\rho \cos\alpha \\
0 \\
\cos\rho \sin\alpha \cos\psi\\
\cos\rho \sin\alpha \sin\psi\\
\end{pmatrix}
\]
Note that the orbit space is homeomorphic to $\mathbb S^3$ and that closed geodesics in the orbit space correspond to minimal tori in $\mathbb S^4$.
Since a complete set of integrals is found for the geodesic flow, the action angle variables allow the study of closed orbits.
Indeed, the geodesic equations admit three independent integrals: $H_0$ which is the velocity squared, $H_1$ and $H_2$ which is the angle between
the geodesic and the $\psi = const$ lines. The subspaces $\psi = const$ (with $H_2=0$) form totally geodesic $2$-spheres in the orbit space.

\begin{prop}
\label{propo2}
The solutions of (\ref{odesystem}) lie on the totally geodesic $2$-sphere defined by $\psi = 0$.
\end{prop}
\begin{proof}
We set the following identifitcations between the 2-sphere $\psi = 0$ and the functions in (\ref{eigen}),
\[
\begin{cases}
\varphi_0 = \sin\rho\\
\varphi_1 = \cos\rho \cos\alpha\\
\varphi_2 = \cos\rho \sin\alpha\\
\end{cases}
\]
The geodesics on this $2$-sphere satisfy our system (\ref{odesystem}), after a change of variable from the geodesic parameter $s$ to $y$ given by,
\[
\frac{ds}{dy} = m^2\varphi_1^2+n^2\varphi_2^2
\]
\end{proof}

The Lawson's bipolar surfaces $\tilde \tau_{r,k}$ are also found in
\cite{15}, where the $(m,n)$ parameters of $R(\theta)$ must be
chosen according to the given pair $(r,k)$. It is shown that they
lie on the ellipse $2\varphi_1^2+\frac{2n^2}{n^2+m^2}\varphi_0^2 =
1$ and that they correspond to an extremal point of $H_1$ (with $H_0
= 1$, $H_2 = 0$). Note that in our situation these surfaces were
also extremals of our first integral $E_1$.

\smallskip
\noindent{\bf Remark.}
A general study of equivariant harmonic maps in spheres is found in \cite{16}. The conformal equivariant immersions of $\mathbb S^1 \times \mathbb S^1$,
given by harmonic maps are shown to correspond to minimal tori.

The system (\ref{odesystem}) is found by considering harmonic immersions in $\mathbb S^4$ and $\mathbb S^1$-equivariance under the transformations $R(\theta)$ above,
after setting to zero the angular momenta $Q_{\alpha\beta}$. The additional
integrals $E_0$, $E_m$ and $E_n$ mentioned in this paper correspond to
our first integral $E_1$ (in our situation the three are dependent, see \cite{12} for the case $m = 1$ and $n = 2$).

\textbf{Acknowledgements.} This research was conducted under the
supervision of Iosif Polterovich and supported by the Undergraduate
Student Research Award of the NSERC. The problem was posed by Prof.
Polterovich and I would like to thank him for many useful
discussions.

I am also grateful to Ahmad El Soufi   for providing an addendum to
\cite{11} with enough details to show proposition \ref{propo1}.

\end{document}